\documentclass[12 pt,a4]{article}
\usepackage{amssymb,amsmath,graphicx,enumerate}
\usepackage{longtable,mathrsfs}
\usepackage{hyperref}
\graphicspath {P:\Dole}

\DeclareMathAlphabet{\mathpzc}{OT1}{pzc}{m}{it}

\numberwithin{equation}{section}

\newtheorem{theorem}{Theorem}[section]
\newtheorem{lemma}{Lemma}[section]

\newtheorem{definition}{Definition}[section]

\newtheorem{prop}{Proposition}[section]

\title{k-Hilfer-Prabhakar Fractional Derivatives and Applications}
\author{\bf S. K. Panchal, \qquad Amol D. Khandagale,\\
\bf Pravinkumar V. Dole\\ Department of mathematics, \\ Dr. Babasaheb Ambedkar Marathwada University,\\ Aurangabad-431004 (M.S.) India.\\
E-mail ID - drpanchalsk@gmail.com\\ kamoldsk@gmail.com\\ pvasudeo.dole@gmail.com  }
\date{\it August 2016.}
\begin{document}
\maketitle{\bf{Abstract}}
\paragraph{}In this paper we define the regularized version of k-Prabhakar fractional derivative, k-Hilfer-Prabhakar fractional derivative, regularized version of k-Hilfer-Prabhakar fractional derivative and find their Laplace and Sumudu transforms. Using these results, the relation between k-Prabhakar fractional derivative and its regularized version involving k-Mittag-Leffler function is obtained. Similarly the relation between k-Hilfer-Prabhakar fractional derivative and its regularized version is also obtained. Further, we find the solutions of some problems in physics in which k-Hilfer-Prabhakar fractional derivative and its regularized version are involved.

\textbf{Keywords:} 
\textit{Fractional calculus; k-Hilfer derivatives; k-Prabhakar derivatives; k-Mittag-Leffler function; Integral transforms.}\\

\textbf{2010 AMS Subject Classification:}\quad $26A33$, $42A38$, $42B10$, $44A10$.
\section{Introduction}
\paragraph{}The fractional calculus is the field of mathematical analysis which deals with the investigation and applications of integrals and derivatives of arbitrary order. In the literature several different definitions of fractional integrals and derivatives are available, like Riemann-Liouville integral, Riemann-Liouville fractional derivative, Caputo fractional derivatives, Liouville integral, Liouville fractional derivative, Miller-Ross sequential fractional derivatives etc. R. Hilfer has defined Hilfer fractional dearivative in \cite{HR}, it is generalization of Riemann-Liouville fractional derivative and Caputo fractional derivative. The Prabhakar integral is obtained from Riemann-Liouville integral operater by extending its kernel involving three parameter Mittag-Leffler function \cite{RRFZ, KSS, PTR}. The Prabhakar derivative is defined by replacing Riemann-Liouville integral operator in  Riemann-Liouville derivative by the Prabhakar integral operator \cite{RRFZ, KSS}. Roberto Garra et. al defined Hilfer-Prabhakar fractional derivative from defination of Hilfer fractional derivative by replacing Reimann Liouvlle integral operator by Prabhakar integral operator and also defined regularized version of Prabhakar and Hilfer-Prabhakar fractional derivatives in \cite{RRFZ}. G. A. Dorrego and R. A. Cerutti defined k-Hilfer fractional derivative in \cite{DC2}. Recently, G. A. Dorrego defined the k-Prabhakar integral and k-Prabhakar fractional derivative in \cite{D}. The kernel of k-Prabhakar integral involve k-Mittag-Leffler function defined by G. A. Dorrego and R. A. Cerutti in \cite{DC1}.\\
In this paper we define regularized version of k-Prabhakar fractional derivative, k-Hilfer-Prabhakar fractional derivative and its regularized version and find their Laplace and Sumudu transforms. The Laplace and Sumudu transforms of k-Prabhakar fractional derivative defined in \cite{D} are obtained. Using these results, the relation between k-Prabhakar fractional derivative and its regularized version, also the relation between k-Hilfer-Prabhakar fractional derivative and its regularized version involving k-Mittag-Leffler function are obtained. The Sumudu transform of k-Prabhakar integral is also obtained.  Lastly, the Laplace transform, Sumudu transform and Fourier transform techniques are used to obtain the solutions of non-homogeneous Cauchy type problems involving the k-Hilfer-Prabhakar fractional derivative and its regularized version. The results obtained in this paper are generalizations of some results obtained in \cite{RRFZ} and \cite{PKD}.
\section{Preliminaries}
\paragraph{} In this section we provide some definitions, theorem, lemma and proposition, which are used in the paper. 

\begin{definition}\label{Def1} 
The Laplace transform of $f(x)$ denoted by $\mathcal{L}[f(x)](s)= F(s)$, is defined as, 
\begin{equation}
\mathcal{L}[f(x)](s)=F(s)=\int_{0}^{\infty} e^{-sx}f(x)dx,\label{eq2.1}
\end{equation}
if the integral in \eqref{eq2.1} exists. The function $f(x)$ is inverse Laplace transform of $F(s)$ and denoted by $\mathcal{L}^{-1}[F(s)](x)=f(x)$.
\end{definition}

\begin{definition}\label{Def2}
The Fourier transform of $f(x)$ is denoted by $\mathcal{F}[f(x)](p)= F(p)$, $p \in \mathbb{R}$, and defined by the integral 
\begin{align}
\mathcal{F}[f(x)](p) = \int_{-\infty}^{\infty} e^{-ipx} f(x) dx.\label{eq2.2}
\end{align}
The inverse Fourier transform of $F(p)$ denoted by $\mathcal{F}^{-1}(F(p))=f(x)$, is defined as,
\begin{align}
f(x) = \mathcal{F}^{-1}[F(p)]= \frac{1}{2\pi}\int_{-\infty}^{\infty} F(p) e^{ipx} dp.\label{eq2.3}
\end{align}
\end{definition}

\begin{definition}\cite{DC1}\textbf{(k-Mittage-Leffler function)}\label{Def3}
Let $ n\in\mathbb{N} $, $ k\in\mathbb{R^{+}} $; $ \alpha, \mu, \gamma\in\mathbb{C} $, $ Re(\alpha)>0 $, $ Re(\mu)>0 $. The k-Mittage-Leffler function is defined as,
\begin{equation}
E_{k, \alpha, \mu}^{\gamma}(z) = \sum_{n=0}^{\infty}\frac{(\gamma)_{n, k}} {\Gamma_{k}(\alpha n+\mu)} \frac{z^{n}}{n!},\label{eq2.4}
\end{equation}
where $(\gamma)_{n, k} =\gamma(\gamma+k)(\gamma+2k)... (\gamma+(n-1)k) $ is the k-Pochhammer symbol and
$ \Gamma_{k}(\alpha)=\int_{0}^{\infty} e^{\frac{-t^{k}}{k}} t^{\alpha-1} dt $ is the k-gamma function \cite{DP}.
\end{definition}

\begin{definition}\cite{DC2}\textbf{(k-Hilfer Derivative)}\label{Def4}\\
Let $ k\in\mathbb{R^{+}} $; $ \mu, \nu\in\mathbb{R} $; $ 0 < \mu < 1 $; $ 0 \leq \nu \leq 1 $. The k-Hilfer Derivatives of $f(x)$ denoted by $ ^{k}D^{\mu, \nu}f(x) $, is defined as,
\begin{equation}
^{k}D^{\mu, \nu}f(x)=\bigg{(}I_{k}^{\nu(\mu-1)}\frac{d}{dt}\big{(}I_{k}^{(1-\nu)(\mu-1)}f\big{)}\bigg{)}(x).\label{eq2.5}
\end{equation}
where,
\begin{equation}
  I_{k}^{\alpha}f(x) = \frac{1}{k\Gamma_{k}(\alpha)}\int_{0}^{x}(x-t)^{\frac{\alpha}{k}-1}f(t)dt \label{eq2.6}
\end{equation}
is k-Riemann-Liouville fractional integral \cite{SK}.
\end{definition}

\begin{definition}\cite{D}{\bf(k-Prabhakar Integral)}\label{Def5}\\
Let $ f \in L^{1} [0, b], 0 < x < b < \infty$ and $k\in\mathbb{R^{+}}$. The k-prabhakar integral is defined as,
\begin{align}
\bigg{(} {_{k}}P_{\alpha, \mu, \omega}^{\gamma} f\bigg{)}(x)= \int_{0}^{x}\frac{(x-t)^{\frac{\mu}{k}-1}}{k} E_{k, \alpha, \mu}^{\gamma}[\omega(x-t)^{\frac{\alpha}{k}}]f(t) dt=( _{k}\varepsilon_{\alpha, \mu, \omega}^{\gamma}* f)(x),\label{eq2.7}
\end{align} 
where * denotes the convolution operation; $\alpha, \mu, \omega, \gamma\,\,\in \,\,\mathbb{C}$; $ Re(\alpha) > 0$, $ Re(\mu)>0 $ and
\begin{align}
_{k}\varepsilon_{\alpha, \mu, \omega}^{\gamma}(t)= \left\{ {\begin{array}{*{20}{l}}{\frac{t^{\frac{\mu}{k}-1}}{k} E_{k, \alpha, \mu}^{\gamma}(\omega t^{\frac{\alpha}{k}}),\qquad\qquad t > 0}\label{eq2.8}\\
{0,\qquad\qquad\qquad\qquad\qquad t \leq 0}.
\end{array}} \right.
\end{align}
For  $\gamma = 0 $, $\bigg{(} {_{k}}P_{\alpha, \mu, \omega}^{0} f\bigg{)}(x)=\bigg{(}I_{k}^{\mu}\phi\bigg{)}(x) $
and for $ \gamma=\mu = 0 $, $\bigg{(} {_{k}}P_{\alpha, 0, \omega}^{0} f\bigg{)}(x)=\phi(x) $. 
\end{definition}

\begin{definition}\cite{D}\textbf{(k-Prabhakar Derivative)}\label{Def6}\\
Let $ k\in\mathbb{R^{+}} $, $\rho, \mu, \omega, \gamma\,\,\in \,\,\mathbb{C};  Re(\alpha), Re(\mu) > 0; m=\big{[}\frac{\mu}{k}\big{]}+ 1$ and $f \in L^{1} [0, b],\\ 0 < t < b < \infty$. The k-Prabhakar derivative of $f(t)$ of order $\mu $ is defined as,
\begin{equation}
_{k}D_{\rho, \mu, \omega}^{\gamma} f(t)= \bigg{(}\frac{d}{dt}\bigg{)}^{m}k^{m}{_{k}}P_{\rho,mk-\mu, \omega}^{-\gamma}f(t).\label{eq2.9}
\end{equation}
\end{definition}

\begin{theorem} Let $F(u)$ and $G(u)$ be Laplace transforms of $f(t)$ and $g(t)$ respectively, then the Laplace transforms of convolution of $f$ and $g$ is
\begin{equation}
\mathcal{L}\big{[}(f*g)(t)\big{]}(u)=F(u)G(u),\label{eq2.10}
\end{equation}
where the convolution is defined as $(f*g)(t)= \int_{0}^{t}f(t-u) g(u) du $.
\end{theorem}

\begin{lemma}\cite{D} Let $\alpha, \mu, \omega, \gamma \in \mathbb{C}$; $Re(\alpha) > 0, Re(\mu) > 0, k\in\mathbb{R^{+}}$ and $\left|\omega k (ku)^{\frac{-\alpha}{k}}\right|<1 $. Then the Laplace transform of the function $_{k}\varepsilon_{\alpha, \mu, \omega}^{\gamma}(t)$ defined in \eqref{eq2.8} is
\begin{equation}
\mathcal{L}\bigg{(}{_{k}\varepsilon}_{\alpha, \mu, \omega}^{\gamma}(t)\bigg{)}(u)= (ku)^{\frac{-\mu}{k}}\bigg{(}1-\omega k (ku)^{\frac{-\alpha}{k}}\bigg{)}^{\frac{-\gamma}{k}}.\label{eq2.11}
\end{equation}
\end{lemma}

\begin{prop}\cite{D}\textbf{(Composition property of k-Prabhakar integral operator)}\label{prop2.1}\\
Let $\alpha, \mu, \nu, \gamma, \sigma, \omega \in \mathbb{C} $; $ k \in \mathbb{R^{+}}$ and $ Re(\alpha)>0, Re(\mu)>0, Re(\nu) > 0 $; then for any $\phi \in L^{1}[0,b]$ and $0 < x < b \leq \infty$ we have,
\begin{align}
\bigg{(}{_{k}}P_{\alpha, \mu, \omega}^{\gamma}{_{k}}P_{\alpha, \nu, \omega}^{\sigma} \phi\bigg{)}(x)= {_{k}}P_{\alpha, \mu+\nu, \omega}^{\gamma+\sigma}\phi(x).\label{eq2.12}
\end{align}
\end{prop}

\begin{definition}\cite{WGK}
Consider a set A defined as,\\
\begin{align}
A=\bigg{\lbrace}f(t)/ \exists\,\, M, \,\, \tau_{1}, \tau_{2} > 0, |f(t)|\leq Me^{\frac{|t|}{\tau_{j}}}\,\, if \,\, t\,\,\in \,\,(-1)^{j}\times[0, \infty)\bigg{\rbrace}\label{2.13}
\end{align}
For all real $t\geq 0$ the Sumudu transform of function $f(t)\in A$ is defined as,
\begin{align}
\mathcal{S}[f(t)](u) = \int_{0}^{\infty}\frac{1}{u}e^{-\frac{t}{u}}f(t) dt,\,\,\, u\in(-\tau_{1}, \tau_{2})\label{2.14}
\end{align}
and is denoted by $F(u) = \mathcal{S}[f(t)](u)$.
\end{definition}

\begin{definition}\cite{WGK}
The function $f(t)$ in \eqref{2.13} is called inverse Sumudu transform of $F(u)$ and is denoted by, 
\begin{align}
f(t) = \mathcal{S}^{-1}[F(u)](t)\label{2.15}
\end{align} 
and the inversion formula for Sumudu transform is given by \cite{WGK} ,
\begin{align}
f(t) = \mathcal{S}^{-1}[F(u)](t)= \frac{1}{2\pi i}\int_{\gamma-i\infty}^{\gamma-i\infty} \frac{1}{u} e^{\frac{t}{u}}f(u) du \label{2.16}
\end{align}
For $Re(\frac{1}{u}) > \gamma$ and $\gamma\,\,\in \,\,\mathbb{C}$.
\end{definition}
\begin{theorem}\cite{BK}
Let $f(t)$ be in $A$ defined in \eqref{2.13}, and let $G_{n}(u)$ denote the Sumudu transform of the $n^{th}$ derivative, $f^{(n)}(t)$ of $f(t)$, then for $n\geq 1$, 
\begin{align}
G_{n}(u)=\frac{G(u)}{u^{n}}-\sum_{k=0}^{n-1}\frac{f^{(n)}(0)}{u^{n-k}}=\frac{G(u)}{u^{n}}-\sum_{k=0}^{n-1}\frac{f^{(n-k-1)}(0)}{u^{k+1}},\label{2.17}
\end{align}
where $G(u)$ denotes the Sumudu transform of $f(t)$.
\end{theorem}
\begin{theorem}\cite{WGK} Let $F(u)$ and $G(u)$ be Sumudu transforms of $f(t)$ and $g(t)$ respectively. The Sumudu transforms of convolution of $f$ and $g$ is
\begin{align}
S\big{[}(f*g)(t)\big{]}(u)=uF(u)G(u),\label{2.18}
\end{align}
where the convolution is defined as $(f*g)(t)=\int_{0}^{t}f(t)g(t-\tau)d\tau$. 
\end{theorem}
\begin{prop}\cite{DC1} Let $\gamma\in \mathbb{C}$ and $k \in \mathbb{R}$. Then the following identity hold,
\begin{align}
\Gamma{_{k}(\gamma)}=k^{\frac{\gamma}{k}-1}\Gamma(\gamma).\label{2.19}
\end{align}
\end{prop}

Now, for $n\in \mathbb{N}$,
\begin{center}
$AC^{n}[a, b] = \bigg{\lbrace} f: [a, b]\rightarrow \mathbb{R}\quad \big{|}\quad \frac{d^{n-1}}{dt^{n-1}}f(t) \in AC[a, b] \bigg{\rbrace}$,
\end{center}
where $AC[a, b]$ denotes the space of all real valued absolutely continuous functions defined on $[a, b]$.

\section{Main Result}In this section, we define the regularized version of k-Prabhakar fractional derivative, k-Hilfer-Prabhakar fractional Derivative and its regularized version and find their Laplace and Sumudu transforms. The Laplace and Sumudu transform of k-Prabhakar fractional derivative are obtained. These results are used to obtain the relation between k-Prabhakar fractional derivative and its regularized version and also the relation between k-Hilfer-Prabhakar fractional derivative and its regularized version involving k-Mittag-Leffler function. The Sumudu transform of k-Prabhakar integral is also obtained. For $k=1$ these results and relations reduces to some results and relations obtained in \cite{RRFZ} and \cite{PKD}.

\begin{definition}\textbf{(Regularized version of k-Prabhakar Fractional Derivative)}\label{Def3.1}\\
Let $ k\in\mathbb{R^{+}} $; $\rho, \mu, \omega, \gamma\,\,\in \,\,\mathbb{C};  Re(\alpha), Re(\mu) > 0; m=\big{[}\frac{\mu}{k}\big{]}+ 1$ and $f \in AC^{m}[0, b],\\ 0 < t < b < \infty$. The regularized version of k-Prabhakar fractional derivative of order $\mu$ of $ f(t) $ denoted by $_{k}^{C}D_{\rho, \mu, \omega}^{\gamma}f(t)$, is defined as,
\begin{equation} 
_{k}^{C}D_{\rho, \mu, \omega}^{\gamma}f(t)=k^{m}{_{k}}P_{\rho,mk-\mu, \omega}^{-\gamma}\bigg{(}\frac{d}{dt}\bigg{)}^{m}f(t).\label{eq3.1}
\end{equation}
\end{definition}

\begin{lemma}\label{lem3.1}
The Laplace transform of k-Prabhakar fractional derivative \eqref{eq2.9} is,
\begin{align*}
&\mathcal{L}\bigg{(}{_{k}D}_{\alpha,\mu, \omega}^{\gamma} f(t)\bigg{)}(u)\\
&=(ku)^{\frac{\mu}{k}}\bigg{(}1-\omega k(ku)^{\frac{-\alpha}{k}}\bigg{)}^{\frac{\gamma}{k}}F(u)-\sum_{n=0}^{m-1}k(ku)^{n}\bigg{[}{_{k}D}_{\alpha,\mu-(n+1)k, \omega}^{\gamma}f(0^{+})\bigg{]}.
\end{align*}
For the case $[\frac{\mu}{k}] + 1 = m = 1$,
\begin{align}
&\mathcal{L}\bigg{(}{_{k}D}_{\alpha,\mu, \omega}^{\gamma} f(t)\bigg{)}(u)\nonumber\\
&=(ku)^{\frac{\mu}{k}}\bigg{(}1-\omega k(ku)^{\frac{-\alpha}{k}}\bigg{)}^{\frac{\gamma}{k}}F(u)-k\bigg{[}{_{k}P}_{\alpha,(k-\mu),\omega}^{-\gamma}f(t)\bigg{]}_{t=0^{+}}.\label{eq3.2}
\end{align}
with $|\omega k(ku)^{\frac{-\alpha}{k}}| < 1$.
\end{lemma}
\textbf{Proof:} Taking Laplace transforms of k-Prabhakar fractional derivative \eqref{eq2.9} and using \eqref{eq2.7}, \eqref{eq2.8}, \eqref{eq2.10} and \eqref{eq2.11}, we get,
\begin{align*}
&\mathcal{L}\bigg{(}{_{k}D}_{\alpha,\mu, \omega}^{\gamma} f(t)\bigg{)}(u)=\mathcal{L}\bigg{(}\frac{d^{m}}{dt^{m}}k^{m}{_{k}P}_{\alpha,(mk-\mu), \omega}^{-\gamma}f(t)\bigg{)}(u)\\
&=k^{m}u^{m}\mathcal{L}\bigg{(}\bigg{(}{_{k}\varepsilon}_{\alpha,(mk-\mu),\omega}^{-\gamma}*f\bigg{)}(t)\bigg{)}(u)\\
&\quad-k^{m}\sum_{n=0}^{m-1}u^{n}\bigg{[}(\frac{d}{dt})^{m-n-1}{_{k}P}_{\alpha,(mk-\mu), \omega}^{-\gamma}f(t)\bigg{]}_{t=o+}\\
&=(ku)^{m}\mathcal{L}\bigg{(}\frac{t^{\frac{(mk-\mu)}{k}-1}}{k} E_{k, \alpha,(mk-\mu)}^{-\gamma}(\omega t^{\frac{\alpha}{k}})\bigg{)}(u)\mathcal{L}\bigg{(}f(t)\bigg{)}(u)\\
&\quad-\sum_{n=0}^{m-1}k^{n+1}u^{n}\bigg{[}k^{m-n-1}\bigg{(}\frac{d}{dt}\bigg{)}^{m-n-1}{_{k}P}_{\alpha,(m-n-1)k-\mu+(n+1)k, \omega}^{-\gamma}f(0^{+})\bigg{]}\\
&=(ku)^{\frac{\mu}{k}}\bigg{(}1-\omega k(ku)^{\frac{-\alpha}{k}}\bigg{)}^{\frac{\gamma}{k}}F(u)-\sum_{n=0}^{m-1}k^{n+1}u^{n}\bigg{[}{_{k}D}_{\alpha,\mu-(n+1)k, \omega}^{\gamma}f(0^{+})\bigg{]}.
\end{align*}
For the case $[\frac{\mu}{k}]+ 1 = m = 1$, we have,
\begin{align*}
&\mathcal{L}\bigg{(}{_{k}D}_{\alpha,\mu, \omega}^{\gamma} f(t)\bigg{)}(u)
&=(ku)^{\frac{\mu}{k}}\bigg{(}1-\omega k(ku)^{\frac{-\alpha}{k}}\bigg{)}^{\frac{\gamma}{k}}F(u)-k\bigg{[}{_{k}P}_{\alpha,(k-\mu), \omega}^{-\gamma}f(t)\bigg{]}_{t=o+}.
\end{align*}
Hence the proof.

\begin{lemma}\label{lem3.2}
The Laplace transform of regularized version of k-Prabhakar fractional derivative \eqref{eq3.1} is,
\begin{align}
\mathcal{L}\bigg{(}{_{k}^{C}D}_{\alpha,\mu, \omega}^{\gamma} f(t)\bigg{)}(u)&=(ku)^{\frac{\mu}{k}}\bigg{(}1-\omega k(ku)^{\frac{-\alpha}{k}}\bigg{)}^{\frac{\gamma}{k}}F(u)\nonumber\\
&\quad-\sum_{n=0}^{m-1}k^{n+1}(ku)^{\frac{-(n+1)k+\mu}{k}}\bigg{(}1-\omega k(ku)^{\frac{-\alpha}{k}}\bigg{)}^{\frac{\gamma}{k}}f^{(n)}(0^{+})\label{eq3.3}
\end{align}
with $|\omega k(ku)^{\frac{-\alpha}{k}}| < 1$.
\end{lemma}
\textbf{Proof:} Taking Laplace transform of regularized version of k-Prabhakar fractional derivative \eqref{eq3.1} and using \eqref{eq2.7}, \eqref{eq2.8}, \eqref{eq2.10}, \eqref{eq2.11}, we get,
\begin{align*}
\mathcal{L}\big{(}{_{k}^{C}D}_{\alpha,\mu, \omega}^{\gamma} f(t)\big{)}(u)&=\mathcal{L}\big{[}k^{m}{_{k}P}_{\alpha,(mk-\mu), \omega}^{-\gamma}\frac{d^{m}}{dt^{m}}f(t)\big{]}(u)\\
&=k^{m}\mathcal{L}\bigg{(}\bigg{(}{_{k}\varepsilon}_{\alpha,(mk-\mu),\omega}^{-\gamma}*\frac{d^{m}f}{dt^{m}}\bigg{)}(t)\bigg{)}(u)\\
&=k^{m}\mathcal{L}\bigg{(}\frac{t^{\frac{(mk-\mu)}{k}-1}}{k} E_{k, \alpha,(mk-\mu)}^{-\gamma}(\omega t^{\frac{\alpha}{k}})\bigg{)}(u)\mathcal{L}\bigg{(}\frac{d^{m}f}{dt^{m}}\bigg{)}(u)\\
&=(ku)^{m}(ku)^{\frac{-(mk-\mu)}{k}}\big{(}1-\omega k(ku)^{\frac{-\alpha}{k}}\big{)}^{\frac{\gamma}{k}}F(u)\\
&-\sum_{n=0}^{m-1}k^{m}u^{m-n-1}(ku)^{\frac{-(mk-\mu)}{k}}\big{(}1-\omega k(ku)^{\frac{-\alpha}{k}}\big{)}^{\frac{\gamma}{k}}f^{(n)}(0^{+}).
\end{align*}
After simplification and rearrangement, we get the required result \eqref{eq3.3}.\\
For absolutely continuous function $f\in AC^{1}[0, b]$,
\begin{align}
\bigg{[}{_{k}P}_{\alpha,(k-\mu), \omega}^{-\gamma}f(t)\bigg{]}_{t=o^{+}}=0\label{3.4} 
\end{align} 
therefore equation \eqref{eq3.2} becomes,
\begin{align}
\mathcal{L}\bigg{(}{_{k}D}_{\alpha,\mu, \omega}^{\gamma} f(t)\bigg{)}(u)=(ku)^{\frac{\mu}{k}}\bigg{(}1-\omega k(ku)^{\frac{-\alpha}{k}}\bigg{)}^{\frac{\gamma}{k}}F(u)\label{eq3.5}
\end{align} 
Thus for $f\in AC^{1}[0, b]$,
\begin{align}
\mathcal{L}\bigg{(}{_{k}^{C}D}_{\alpha,\mu, \omega}^{\gamma} f(t)\bigg{)}(u)&=\mathcal{L}\bigg{(}{_{k}D}_{\alpha,\mu, \omega}^{\gamma} f(t)\bigg{)}(u)\nonumber\\
&\quad-\sum_{n=0}^{m-1}k^{n+1}(ku)^{\frac{-(n+1)k+\mu}{k}}\bigg{(}1-\omega k(ku)^{\frac{-\alpha}{k}}\bigg{)}^{\frac{\gamma}{k}}f^{(n)}(0^{+}).\label{eq3.6}
\end{align}
Taking inverse Laplace transform of \eqref{eq3.6}, we get, the relation between k-Prabhakar fractional derivative and its regularized version in terms of k-Mittag-Leffler function as below,
\begin{align}
{_{k}^{C}D}_{\alpha,\mu, \omega}^{\gamma} f(t)={_{k}D}_{\alpha,\mu, \omega}^{\gamma} f(t)-\sum_{n=0}^{m-1}k^{n}t^{\frac{nk-\mu}{k}}E_{k, \alpha,(n+1)k-\mu}^{-\gamma}(\omega t^{\frac{\alpha}{k}})f^{(n)}(0^{+}),\label{eq3.7}
\end{align}
for $f\in AC^{1}[0, b]$.

\begin{definition}{\bf(k-Hilfer-Prabhakar Fractional Derivative)}\label{Def3.2}\\
Let $f \in L^{1}[a, b], 0 < t < b < \infty$; $ k, \alpha > 0; \gamma,\omega\in \mathbb{R}; \mu\in(0, 1),\nu\in [0, 1]$ and $(f* {_{k}}\varepsilon_{\rho,(1-\nu)(k-\mu),\omega}^{-\gamma (1-\nu)})(t) \in AC^{1}[0, b]$. The k-Hilfer-Prabhakar fractional derivative of $f(t)$ of order $\mu$ denoted by $_{k}\mathcal{D}_{\alpha, \omega, 0^{+}}^{\gamma, \mu, \nu} f(t)$ is defined as,
\begin{align}
_{k}\mathcal{D}_{\alpha, \omega, 0^{+}}^{\gamma, \mu, \nu} f(t)= k\bigg{(} {_{k}}P_{\alpha, \nu(k-\mu), \omega, 0^{+}}^{-\gamma\nu}\frac{d}{dt}\big{(} {_{k}}P_{\alpha,(1-\nu)(k-\mu), \omega, 0^{+}}^{-\gamma(1-\nu)}f\big{)}\bigg{)}(t).\label{eq3.8}  
\end{align}
\end{definition}

We observe that for $ \gamma = 0 $, the k-Hilfer-Prabhakar fractional derivative \eqref{eq3.8}  reduces to k-Hilfer fractional derivative \eqref{eq2.5} of order $\mu$. For $ \nu = 0 $ and $ \nu = 1 $ it reduces to k-Prabhakar fractional derivative \eqref{eq2.9} and regularized version of k-Prabhakar fractional derivative \eqref{eq3.1} of order $\mu$ respectively (note $m=1$).\\
	In order to consider the Cauchy problems in which initial conditions only depends on the function and its integer-order derivatives, we use the regularized version of k-Hilfer-Prabhakar fractional derivative defined as below.

\begin{definition}{\bf(Regularized Version of k-Hilfer-Prabhakar Fractional Derivative)}\label{Def3.3}\\
Let $f \in AC^{1}[0, b], 0 < t < b < \infty$ and $ k, \alpha > 0; \gamma,\omega\in \mathbb{R}; \mu\in(0, 1),\nu\in [0, 1]$. The regularized version of k-Hilfer-Prabhakar fractional derivative of $f(t)$ denoted by $_{k}^{C}\mathcal{D}_{\alpha,\omega,0^{+}}^{\gamma,\mu,\nu}f(t)$ is defined as,
\begin{align}
_{k}^{C}\mathcal{D}_{\alpha,\omega,0^{+}}^{\gamma,\mu,\nu}f(t)=k\bigg{(}{_{k}P}_{\alpha, \nu(k-\mu), \omega, 0^{+}}^{-\gamma\nu}{_{k}P}_{\alpha,(1-\nu)(k-\mu), \omega, 0^{+}}^{-\gamma(1-\nu)}\frac{d}{dt}f\bigg{)}(t),\label{eq3.9}
\end{align}
Using the Proposition \eqref{prop2.1}, the regularized version of k-Hilfer-Prabhakar fractional derivative can be expressed as,
\begin{align}
_{k}^{C}\mathcal{D}_{\alpha,\omega,0^{+}}^{\gamma,\mu,\nu}f(t)=k\bigg{(}{_{k}P}_{\alpha, (k-\mu), \omega, 0^{+}}^{-\gamma}\frac{d}{dt}f\bigg{)}(t).\label{eq3.10}
\end{align}
\end{definition}

\begin{lemma}\label{lem3.3}
The Laplace transform of k-Hilfer-Prabhakar fractional derivative \eqref{eq3.8} is,
\begin{align}
\mathcal{L}\bigg{(}{_{k}\mathcal{D}}_{\alpha, \omega, 0^{+}}^{\gamma, \mu, \nu} f(t)\bigg{)}(u)&=(ku)^{\frac{\mu}{k}}\bigg{(}1-\omega k(ku)^{\frac{-\alpha}{k}}\bigg{)}^{\frac{\gamma}{k}}F(u)-k(ku)^{\frac{-\nu(k-\mu)}{k}}\nonumber\\
&\qquad\bigg{(}1-\omega k(ku)^{\frac{-\alpha}{k}}\bigg{)}^{\frac{\gamma\nu}{k}}\bigg{[}{_{k}P}_{\alpha,(1-\nu)(k-\mu),\omega, 0^{+}}^{-\gamma(1-\nu)}f(t)\bigg{]}_{t=0^{+}},\label{eq3.11}
\end{align}
where $F(u)$ denotes the Laplace transform of $f(t)$.
\end{lemma}
\textbf{Proof:} Taking Laplace transform of k-Hilfer-Prabhakar fractional derivative \eqref{eq3.8} and using \eqref{eq2.7}, \eqref{eq2.8}, \eqref{eq2.10} and \eqref{eq2.11}, we have,
\begin{align*}
&\mathcal{L}\bigg{(}{_{k}\mathcal{D}}_{\alpha, \omega, 0^{+}}^{\gamma, \mu, \nu} f(t)\bigg{)}(u)=\mathcal{L}\bigg{[}k\bigg{(}{_{k}P}_{\alpha,\nu(k-\mu),\omega,0^{+}}^{-\gamma\nu}\frac{d}{dt}\big{(}{_{k}P}_{\alpha,(1-\nu)(k-\mu),\omega,0^{+}}^{-\gamma(1-\nu)}f\big{)}\bigg{)}(t)\bigg{]}(u)\\
&=k\mathcal{L}\bigg{[} \bigg{(}{_{k}\varepsilon}_{\alpha,\nu(k-\mu),\omega}^{-\gamma\nu}*\frac{d}{dt}\big{(}{_{k}P}_{\alpha,(1-\nu)(k-\mu),\omega,0^{+}}^{-\gamma(1-\nu)}f\big{)}\bigg{)}(t)\bigg{]}(u)\\
&=k(ku)^{\frac{-\nu(k-\mu)}{k}}\bigg{(}1-\omega k(ku)^{\frac{-\alpha}{k}}\bigg{)}^{\frac{\gamma\nu}{k}}\\
&\quad\bigg{\lbrace}u\mathcal{L}\big{[}{_{k}P}_{\alpha,(1-\nu)(k-\mu), \omega, 0^{+}}^{-\gamma(1-\nu)}f(t)\big{]}(u)-\big{[}{_{k}P}_{\alpha,(1-\nu)(k-\mu), \omega, 0^{+}}^{-\gamma(1-\nu)}f(t)\big{]}_{t=0^{+}}\bigg{\rbrace}\\
&=ku(ku)^{\frac{-\nu(k-\mu)}{k}}\bigg{(}1-\omega k(ku)^{\frac{-\alpha}{k}}\bigg{)}^{\frac{\gamma\nu}{k}}\mathcal{L}\bigg{(} \big{(}{_{k}\varepsilon}_{\alpha,(1-\nu)(k-\mu), \omega}^{-\gamma(1-\nu)}*f\big{)}(t)\bigg{)}(u)\\
&\quad-k(ku)^{\frac{-\nu(k-\mu)}{k}}\bigg{(}1-\omega k(ku)^{\frac{-\alpha}{k}}\bigg{)}^{\frac{\gamma\nu}{k}}\bigg{[}{_{k}P}_{\alpha,(1-\nu)(k-\mu), \omega, 0^{+}}^{-\gamma(1-\nu)}f(t)\bigg{]}_{t=0^{+}}\\
&=(ku)^{\frac{(1-\nu)k+\nu\mu}{k}}\bigg{(}1-\omega k(ku)^{\frac{-\alpha}{k}}\bigg{)}^{\frac{\gamma\nu}{k}}(ku)^{\frac{-(1-\nu)(k-\mu)}{k}}\bigg{(}1-\omega k(ku)^{\frac{-\alpha}{k}}\bigg{)}^{\frac{\gamma(1-\nu)}{k}}F(u)\\
&\quad-k(ku)^{\frac{-\nu(k-\mu)}{k}}\bigg{(}1-\omega k(ku)^{\frac{-\alpha}{k}}\bigg{)}^{\frac{\gamma\nu}{k}}\bigg{[}{_{k}P}_{\alpha,(1-\nu)(k-\mu),\omega, 0^{+}}^{-\gamma(1-\nu)}f(t)\bigg{]}_{t=0^{+}}
\end{align*}
After simplification, we get \eqref{eq3.11} as desired.
\begin{lemma}\label{lem3.4}
The Laplace transforms of the regularized version of k-Hilfer-Prabhakar fractional derivative\eqref{eq3.9} of order $\mu$ is,
\begin{align}
&\mathcal{L}\bigg{(}{ _{k}^{C}\mathcal{D}}_{\alpha, \omega, 0^{+}}^{\gamma, \mu, \nu} f(t)\bigg{)}(u)\nonumber\\
&=(ku)^{\frac{\mu}{k}}\bigg{(}1-\omega k(ku)^{\frac{-\alpha}{k}}\bigg{)}^{\frac{\gamma}{k}}F(u)-k(ku)^{\frac{-(k-\mu)}{k}}\bigg{(}1-\omega k(ku)^{\frac{-\alpha}{k}}\bigg{)}^{\frac{\gamma}{k}}f(0^{+}),\label{eq3.12}
\end{align}
where $F(u)$ denotes the Laplace transform of $f(t)$.
\end{lemma}
\textbf{Proof:} Taking Laplace transforms of regularized version of k-Hilfer-Prabhakar fractional derivative \eqref{eq3.9} and using \eqref{eq2.7}, \eqref{eq2.8}, \eqref{eq2.10} and \eqref{eq2.11}, we have,\\
\begin{align*}
&\mathcal{L}\bigg{(}{_{k}^{C}\mathcal{D}}_{\alpha, \omega, 0^{+}}^{\gamma, \mu,\nu} f(t)\bigg{)}(u)=\mathcal{L}\bigg{(}k\bigg{(}{_{k}P}_{\alpha,\nu(k-\mu),\omega,0^{+}}^{-\gamma\nu}\big{(}{_{k}P}_{\alpha,(1-\nu)(k-\mu),\omega,0^{+}}^{-\gamma(1-\nu)}\frac{d}{dt}f\big{)}\bigg{)}(t)\bigg{)}(u)\\
&=k\mathcal{L}\bigg{(} \bigg{(}{_{k}\varepsilon}_{\alpha,\nu(k-\mu),\omega}^{-\gamma\nu}*\big{(}{_{k}P}_{\alpha,(1-\nu)(k-\mu),\omega,0^{+}}^{-\gamma(1-\nu)}\frac{d}{dt}f\big{)}\bigg{)}(t)\bigg{)}(u)\\
&=k(ku)^{\frac{-\nu(k-\mu)}{k}}\bigg{(}1-\omega k(ku)^{\frac{-\alpha}{k}}\bigg{)}^{\frac{\gamma\nu}{k}}\mathcal{L}\bigg{[}\bigg{(}{_{k}\varepsilon}_{\alpha,(1-\nu)(k-\mu),\omega}^{-\gamma(1-\nu)}*\frac{d}{dt}f\bigg{)}(t)\bigg{]}(u)\\
&=k(ku)^{\frac{-\nu(k-\mu)}{k}}\big{(}1-\omega k(ku)^{\frac{-\alpha}{k}}\big{)}^{\frac{\gamma\nu}{k}}(ku)^{\frac{-(1-\nu)(k-\mu)}{k}}\big{(}1-\omega k(ku)^{\frac{-\alpha}{k}}\big{)}^{\frac{\gamma(1-\nu)}{k}}\mathcal{L}\bigg{[}\frac{df}{dt}\bigg{]}(u)\\
&=(ku)^{\frac{\mu}{k}}\big{(}1-\omega k(ku)^{\frac{-\alpha}{k}}\big{)}^{\frac{\gamma}{k}}F(u)-k(ku)^{\frac{-(k-\mu)}{k}}\big{(}1-\omega k(ku)^{\frac{-\alpha}{k}}\big{)}^{\frac{\gamma}{k}}f(0^{+})
\end{align*}
Hence the proof.\\
\textbf{Alternating Proof of above Lemma :} Taking Laplace transform of regularized version of k-Hilfer-Prabhakar fractional derivative \eqref{eq3.10} and using \eqref{eq2.7}, \eqref{eq2.8}, \eqref{eq2.10} and \eqref{eq2.11}, we have,\\
\begin{align*}
&\mathcal{L}\bigg{(}{_{k}^{C}\mathcal{D}}_{\alpha, \omega, 0^{+}}^{\gamma,\mu,\nu} f(t)\bigg{)}(u)
=k\mathcal{L}\bigg{(}\bigg{(}{_{k}\varepsilon}_{\alpha,(k-\mu),\omega}^{-\gamma}*\frac{d}{dt}f\bigg{)}(t)\bigg{)}(u)\\
&=k(ku)^{\frac{-(k-\mu)}{k}}\bigg{(}1-\omega k(ku)^{\frac{-\alpha}{k}}\bigg{)}^{\frac{\gamma}{k}}\bigg{(}uF(u)-f(0^{+})\bigg{)}\\
&=(ku)^{\frac{\mu}{k}}\bigg{(}1-\omega k(ku)^{\frac{-\alpha}{k}}\bigg{)}^{\frac{\gamma}{k}}F(u)-k(ku)^{\frac{-(k-\mu)}{k}}\bigg{(}1-\omega k(ku)^{\frac{-\alpha}{k}}\bigg{)}^{\frac{\gamma}{k}}f(0^{+}).
\end{align*}
For absolutely continuous function $f\in AC^{1}[0, b]$,
\begin{align}
\big{[}{_{k}P}_{\rho,(1-\nu)(k-\mu),\omega, 0^{+}}^{-\gamma(1-\nu)}f(t)\big{]}_{t=0^{+}}=0\label{3.13} 
\end{align} 
therefore equation \eqref{eq3.11} becomes,
\begin{align}
\mathcal{L}\bigg{(}{_{k}\mathcal{D}}_{\alpha, \omega, 0^{+}}^{\gamma, \mu, \nu} f(t)\bigg{)}(u)=(ku)^{\frac{\mu}{k}}\bigg{(}1-\omega k(ku)^{\frac{-\alpha}{k}}\bigg{)}^{\frac{\gamma}{k}}F(u).\label{eq3.14}
\end{align} 
Thus, for $f\in AC^{1}[0, b]$, we have,
\begin{align}
\mathcal{L}\bigg{(}{_{k}^{C}\mathcal{D}}_{\alpha, \omega, 0^{+}}^{\gamma,\mu,\nu} f(t)\bigg{)}(u)&=\mathcal{L}\bigg{(}{_{k}\mathcal{D}}_{\alpha, \omega, 0^{+}}^{\gamma, \mu, \nu} f(t)\bigg{)}(u)\nonumber\\
&\quad-k(ku)^{\frac{-(k-\mu)}{k}}\bigg{(}1-\omega k(ku)^{\frac{-\alpha}{k}}\bigg{)}^{\frac{\gamma}{k}}f(0^{+})\label{eq3.15}
\end{align}
Taking inverse Laplace transform of \eqref{eq3.15}, we get the relation between k-Hilfer-Prabhakar fractional derivative and regularized version of k-Hilfer-Prabhakar fractional derivative in terms of k-Mittag-Leffler function as below,
\begin{align}
{_{k}^{C}\mathcal{D}}_{\alpha, \omega, 0^{+}}^{\gamma,\mu,\nu} f(t)={_{k}\mathcal{D}}_{\alpha, \omega, 0^{+}}^{\gamma, \mu, \nu} f(t)-\
t^{\frac{-\mu}{k}} E_{k, \alpha,(k-\mu)}^{-\gamma}(\omega t^{\frac{\alpha}{k}})f(0^{+})\label{eq3.16}
\end{align}
for $f\in AC^{1}[0, b]$.

\begin{lemma}\label{lem3.5}
The Sumudu transform of k-Prabhakar integral \eqref{eq2.7} is
\begin{align}
S\bigg{(} {_{k}}P_{\alpha, \mu, \omega}^{\gamma} f(x)\bigg{)}(u)=u^{-1}{\bigg{(}\frac{u}{k}\bigg{)}^{\frac{\mu}{k}}}\Bigg{(}1-\omega k\bigg{(}\frac{u}{k}\bigg{)}^{\frac{\alpha}{k}}\Bigg{)}^{\frac{-\gamma}{k}} F(u),\label{eq3.17}
\end{align}
provided  $|\omega k\big{(}\frac{u}{k}\big{)}^{\frac{\alpha}{k}}| < 1$, where $F(u)$ is Sumudu transform of $f(x)$.
\end{lemma}
\textbf{Proof:} Taking the Sumudu transform of $_{k}\varepsilon_{\alpha, \mu, \omega}^{\gamma}(x)$ defined in \eqref{eq2.8}, and using \eqref{eq2.4}, \eqref{2.14}, \eqref{2.19}, we have
\begin{align}
S[_{k}\varepsilon_{\alpha, \mu, \omega}^{\gamma}(x)](u)&=\frac{1}{uk}\int_{0}^{\infty}e^{\frac{-x}{u}} x^{\frac{\mu}{k}-1} E_{k, \alpha, \mu}^{\gamma}(\omega x^{\frac{\alpha}{k}})dx \nonumber\\
&=\sum_{n=0}^{\infty}\frac{(\gamma)_{n, k}} {\Gamma_{k}(\alpha n+\mu)} \frac{\omega^{n}}{n!} \frac{1}{uk}\int_{0}^{\infty}e^{\frac{-x}{u}} x^{\frac{n\alpha+\mu}{k}-1}dx \nonumber\\ 
&=\frac{1}{k}\sum_{n=0}^{\infty}\frac{\omega^{n}(\gamma)_{n, k}} {n!k^{\frac{n\alpha+\mu}{k}-1}\Gamma{(\frac{n\alpha+\mu}{k})}}\Gamma{(\frac{n\alpha+\mu}{k})} u^{\frac{n\alpha+\mu}{k}-1} \nonumber\\
&=u^{-1}\bigg{(}\frac{u}{k}\bigg{)}^{\frac{\mu}{k}}\sum_{n=0}^{\infty}\frac{(\gamma)_{n, k}} {n!}\Bigg{(}\omega \bigg{(}\frac{u}{k}\bigg{)}^{\frac{\alpha}{k}}\Bigg{)}^{n} \nonumber\\
&=u^{-1}\bigg{(}\frac{u}{k}\bigg{)}^{\frac{\mu}{k}}\Bigg{(}1-\omega k\bigg{(}\frac{u}{k}\bigg{)}^{\frac{\alpha}{k}}\Bigg{)}^{\frac{-\gamma}{k}},\label{eq3.18}
\end{align}
for $|\omega k\big{(}\frac{u}{k}\big{)}^{\frac{\alpha}{k}}| < 1$.
Now taking Sumudu transform of k-Prabhakar integral \eqref{eq2.7} and using \eqref{2.18}, \eqref{eq3.18}, we get required result \eqref{eq3.17}.
\begin{lemma}\label{lem3.6}
The Sumudu transform of k-Prabhakar fractional derivative \eqref{eq2.9} is,
\begin{align*}
&\mathcal{S}\bigg{(}{_{k}D}_{\alpha,\mu, \omega}^{\gamma} f(t)\bigg{)}(u)\\
&=\bigg{(}\frac{u}{k}\bigg{)}^{\frac{-\mu}{k}}\bigg{(}1-\omega k\bigg{(}\frac{u}{k}\bigg{)}^{\frac{\alpha}{k}}\bigg{)}^{\frac{\gamma}{k}}F(u)-\sum_{n=0}^{m-1}\bigg{(}\frac{k}{u}\bigg{)}^{n+1}\bigg{[}{_{k}D}_{\alpha,\mu-(n+1)k, \omega}^{\gamma}f(0^{+})\bigg{]}.
\end{align*}
For the case $[\frac{\mu}{k}] + 1 = m = 1$,
\begin{align}
&\mathcal{S}\bigg{(}{_{k}D}_{\alpha,\mu, \omega}^{\gamma} f(t)\bigg{)}(u)\nonumber\\
&=\bigg{(}\frac{u}{k}\bigg{)}^{\frac{-\mu}{k}}\bigg{(}1-\omega k\bigg{(}\frac{u}{k}\bigg{)}^{\frac{\alpha}{k}}\bigg{)}^{\frac{\gamma}{k}}F(u)-\frac{k}{u}\bigg{[}{_{k}P}_{\alpha,(k-\mu),\omega}^{-\gamma}f(t)\bigg{]}_{t=0^{+}},\label{eq3.19}
\end{align}
with $|\omega k\big{(}\frac{u}{k}\big{)}^{\frac{\alpha}{k}}| < 1$\\
where $F(u)$ is the Sumudu transform of $f(t)$.
\end{lemma}
\textbf{Proof:} Taking Sumudu transforms of k-Prabhakar fractional derivative \eqref{eq2.9} and using \eqref{2.17},  \eqref{2.18}, \eqref{eq3.18}, we get,
\begin{align*}
&\mathcal{S}\bigg{(}{_{k}D}_{\alpha,\mu, \omega}^{\gamma} f(t)\bigg{)}(u)=\mathcal{S}\bigg{(}\frac{d^{m}}{dt^{m}}k^{m}{_{k}P}_{\alpha,(mk-\mu), \omega}^{-\gamma}f(t)\bigg{)}(u)\\
&=\bigg{(}\frac{k}{u}\bigg{)}^{m}\mathcal{S}\bigg{(}\bigg{(}{_{k}\varepsilon}_{\alpha,(mk-\mu),\omega}^{-\gamma}*f\bigg{)}(t)\bigg{)}(u)\\
&\quad-k^{m}\sum_{n=0}^{m-1}\frac{1}{u^{n+1}}\bigg{[}(\frac{d}{dt})^{m-n-1}{_{k}P}_{\alpha,(mk-\mu), \omega}^{-\gamma}f(t)\bigg{]}_{t=o+}\\
&=\bigg{(}\frac{k}{u}\bigg{)}^{m}uu^{-1}\bigg{(}\frac{u}{k}\bigg{)}^{\frac{mk-\mu}{k}}\bigg{(}1-\omega k\bigg{(}\frac{k}{u}\bigg{)}^{\frac{\alpha}{k}}\bigg{)}^{\frac{\gamma}{k}}\mathcal{S}\bigg{(}f(t)\bigg{)}(u)\\
&\quad-\sum_{n=0}^{m-1}\bigg{(}\frac{k}{u}\bigg{)}^{n+1}\bigg{[}k^{m-n-1}\bigg{(}\frac{d}{dt}\bigg{)}^{m-n-1}{_{k}P}_{\alpha,(m-n-1)k-\mu+(n+1)k, \omega}^{-\gamma}f(0^{+})\bigg{]}\\
&=\bigg{(}\frac{u}{k}\bigg{)}^{\frac{-\mu}{k}}\bigg{(}1-\omega k\bigg{(}\frac{k}{u}\bigg{)}^{\frac{\alpha}{k}}\bigg{)}^{\frac{\gamma}{k}}F(u)-\sum_{n=0}^{m-1}\bigg{(}\frac{k}{u}\bigg{)}^{n+1}\bigg{[}{_{k}D}_{\alpha,\mu-(n+1)k, \omega}^{\gamma}f(0^{+})\bigg{]}.
\end{align*}
For the case $[\frac{\mu}{k}]+ 1 = m = 1$, we have,
\begin{align*}
&\mathcal{S}\bigg{(}{_{k}D}_{\alpha,\mu, \omega}^{\gamma} f(t)\bigg{)}(u)
=\bigg{(}\frac{k}{u}\bigg{)}^{\frac{-\mu}{k}}\bigg{(}1-\omega k(ku)^{\frac{-\alpha}{k}}\bigg{)}^{\frac{\gamma}{k}}F(u)-\frac{k}{u}\bigg{[}{_{k}P}_{\alpha,(k-\mu), \omega}^{-\gamma}f(t)\bigg{]}_{t=o+}.
\end{align*}
Hence the proof.
\begin{lemma}\label{lem3.7}
The Sumudu transform of regularized version of k-Prabhakar fractional derivative \eqref{eq3.1} is,
\begin{align}
\mathcal{S}\bigg{(}{_{k}^{C}D}_{\alpha,\mu, \omega}^{\gamma} f(t)\bigg{)}(u)&=\bigg{(}\frac{u}{k}\bigg{)}^{\frac{-\mu}{k}}\bigg{(}1-\omega k\bigg{(}\frac{k}{u}\bigg{)}^{\frac{\alpha}{k}}\bigg{)}^{\frac{\gamma}{k}}F(u)\nonumber\\
&-\sum_{n=0}^{m-1}k^{n}\bigg{(}\frac{u}{k}\bigg{)}^{\frac{nk-\mu}{k}}\bigg{(}1-\omega k\bigg{(}\frac{k}{u}\bigg{)}^{\frac{\alpha}{k}}\bigg{)}^{\frac{\gamma}{k}}f^{(n)}(0).\label{eq3.20}
\end{align}
with $|\omega k\big{(}\frac{u}{k}\big{)}^{\frac{\alpha}{k}}| < 1$\\
where $F(u)$ is the Sumudu transform of $g(t)$.
\end{lemma}
\textbf{Proof:} Taking Laplace transform of regularized version of k-Prabhakar fractional derivative \eqref{eq3.1} and using \eqref{eq2.7}, \eqref{2.17}, \eqref{2.18}, \eqref{eq3.18}, we get,
\begin{align*}
\mathcal{S}\big{(}{_{k}^{C}D}_{\alpha,\mu,\omega}^{\gamma} f(t)\big{)}(u)
&=k^{m}\mathcal{S}\bigg{(}\bigg{(}{_{k}\varepsilon}_{\alpha,(mk-\mu),\omega}^{-\gamma}*\frac{d^{m}f}{dt^{m}}\bigg{)}(t)\bigg{)}(u)\\
&=k^{m}\bigg{(}\frac{u}{k}\bigg{)}^{\frac{mk-\mu}{k}}\bigg{(}1-\omega k\bigg{(}\frac{k}{u}\bigg{)}^{\frac{\alpha}{k}}\bigg{)}^{\frac{\gamma}{k}}\bigg{[}\frac{F(u)}{u^{m}}-\sum_{n=0}^{m-1}\frac{f^{(n)}(0)}{u^{m-n}}\bigg{]}\\
&=\bigg{(}\frac{u}{k}\bigg{)}^{\frac{-\mu}{k}}\bigg{(}1-\omega k\bigg{(}\frac{k}{u}\bigg{)}^{\frac{\alpha}{k}}\bigg{)}^{\frac{\gamma}{k}}F(u)\\
&-\sum_{n=0}^{m-1}k^{n}\bigg{(}\frac{u}{k}\bigg{)}^{\frac{nk-\mu}{k}}\bigg{(}1-\omega k\bigg{(}\frac{k}{u}\bigg{)}^{\frac{\alpha}{k}}\bigg{)}^{\frac{\gamma}{k}}f^{(n)}(0).
\end{align*}
Hence proved.

\begin{lemma}\label{lem3.8}
The Sumudu transform of k-Hilfer-Prabhakar fractional derivative \eqref{eq3.8} is,
\begin{align}
\mathcal{S}\bigg{(}{_{k}\mathcal{D}}_{\alpha, \omega, 0^{+}}^{\gamma, \mu, \nu} f(t)\bigg{)}(u)&=\bigg{(}\frac{u}{k}\bigg{)}^{\frac{-\mu}{k}}\bigg{(}1-\omega k\bigg{(}\frac{u}{k}\bigg{)}^{\frac{\alpha}{k}}\bigg{)}^{\frac{\gamma}{k}}F(u)-\bigg{(}\frac{u}{k}\bigg{)}^{\frac{\nu(k-\mu)}{k}-1}\nonumber\\
&\qquad\bigg{(}1-\omega k\bigg{(}\frac{u}{k}\bigg{)}^{\frac{\alpha}{k}}\bigg{)}^{\frac{\gamma\nu}{k}}\bigg{[}{_{k}P}_{\alpha,(1-\nu)(k-\mu),\omega, 0^{+}}^{-\gamma(1-\nu)}f(t)\bigg{]}_{t=0^{+}},\label{eq3.21}
\end{align}
where $F(u)$ denotes the Sumudu transform of $f(t)$.
\end{lemma}
\textbf{Proof:} Taking Laplace transform of k-Hilfer-Prabhakar fractional derivative \eqref{eq3.8} and using \eqref{eq2.7}, \eqref{eq2.8}, \eqref{eq2.10} and \eqref{eq2.11}, we have,
\begin{align*}
&\mathcal{S}\bigg{(}{_{k}\mathcal{D}}_{\alpha, \omega, 0^{+}}^{\gamma, \mu, \nu} f(t)\bigg{)}(u)=\mathcal{S}\bigg{[}k\bigg{(}{_{k}P}_{\alpha,\nu(k-\mu),\omega,0^{+}}^{-\gamma\nu}\frac{d}{dt}\big{(}{_{k}P}_{\alpha,(1-\nu)(k-\mu),\omega,0^{+}}^{-\gamma(1-\nu)}f\big{)}\bigg{)}(t)\bigg{]}(u)\\
&=k\mathcal{S}\bigg{[} \bigg{(}{_{k}\varepsilon}_{\alpha,\nu(k-\mu),\omega}^{-\gamma\nu}*\frac{d}{dt}\big{(}{_{k}P}_{\alpha,(1-\nu)(k-\mu),\omega,0^{+}}^{-\gamma(1-\nu)}f\big{)}\bigg{)}(t)\bigg{]}(u)\\
&=k\bigg{(}\frac{u}{k}\bigg{)}^{\frac{\nu(k-\mu)}{k}}\bigg{(}1-\omega k\bigg{(}\frac{u}{k}\bigg{)}^{\frac{\alpha}{k}}\bigg{)}^{\frac{\gamma\nu}{k}}\\
&\quad\frac{1}{u}\bigg{\lbrace}\mathcal{S}\big{[}{_{k}P}_{\alpha,(1-\nu)(k-\mu), \omega, 0^{+}}^{-\gamma(1-\nu)}f(t)\big{]}(u)-\big{[}{_{k}P}_{\alpha,(1-\nu)(k-\mu), \omega, 0^{+}}^{-\gamma(1-\nu)}f(t)\big{]}_{t=0^{+}}\bigg{\rbrace}\\
&=\frac{k}{u}\bigg{(}\frac{u}{k}\bigg{)}^{\frac{\nu(k-\mu)}{k}}\bigg{(}1-\omega k\bigg{(}\frac{u}{k}\bigg{)}^{\frac{\alpha}{k}}\bigg{)}^{\frac{\gamma\nu}{k}}\mathcal{L}\bigg{(} \big{(}{_{k}\varepsilon}_{\alpha,(1-\nu)(k-\mu), \omega}^{-\gamma(1-\nu)}*f\big{)}(t)\bigg{)}(u)\\
&\quad-\frac{k}{u}\bigg{(}\frac{u}{k}\bigg{)}^{\frac{\nu(k-\mu)}{k}}\bigg{(}1-\omega k\bigg{(}\frac{u}{k}\bigg{)}^{\frac{\alpha}{k}}\bigg{)}^{\frac{\gamma\nu}{k}}\bigg{[}{_{k}P}_{\alpha,(1-\nu)(k-\mu), \omega, 0^{+}}^{-\gamma(1-\nu)}f(t)\bigg{]}_{t=0^{+}}\\
&=\frac{k}{u}\bigg{(}\frac{u}{k}\bigg{)}^{\frac{\nu(k-\mu)}{k}}\bigg{(}1-\omega k\bigg{(}\frac{u}{k}\bigg{)}^{\frac{\alpha}{k}}\bigg{)}^{\frac{\gamma\nu}{k}}\bigg{(}\frac{u}{k}\bigg{)}^{\frac{(1-\nu)(k-\mu)}{k}}\bigg{(}1-\omega k\bigg{(}\frac{u}{k}\bigg{)}^{\frac{\alpha}{k}}\bigg{)}^{\frac{\gamma(1-\nu)}{k}}F(u)\\
&\quad-\frac{k}{u}\bigg{(}\frac{u}{k}\bigg{)}^{\frac{\nu(k-\mu)}{k}}\bigg{(}1-\omega k\bigg{(}\frac{u}{k}\bigg{)}^{\frac{\alpha}{k}}\bigg{)}^{\frac{\gamma\nu}{k}}\bigg{[}{_{k}P}_{\alpha,(1-\nu)(k-\mu), \omega, 0^{+}}^{-\gamma(1-\nu)}f(t)\bigg{]}_{t=0^{+}}
\end{align*}
After simplification, we get \eqref{eq3.11} as desired.
\begin{lemma}\label{lem3.9}
The Sumudu transforms of the regularized version of k-Hilfer-Prabhakar fractional derivative\eqref{eq3.9} of order $\mu$ is,
\begin{align}
&\mathcal{S}\bigg{(}{ _{k}^{C}\mathcal{D}}_{\alpha, \omega, 0^{+}}^{\gamma, \mu, \nu} f(t)\bigg{)}(u)=\bigg{(}\frac{u}{k}\bigg{)}^{\frac{-\mu}{k}}\Bigg{(}1-\omega k\bigg{(}\frac{u}{k}\bigg{)}^{\frac{\alpha}{k}}\Bigg{)}^{\frac{\gamma}{k}}\bigg{(}F(u)-f(0^{+})\bigg{)},\label{eq3.22}
\end{align}
where $F(u)$ denotes the Sumudu transform of $f(t)$.
\end{lemma}
\textbf{Proof:}Taking Sumudu transform of regularized version of k-Hilfer-Prabhakar fractional derivative \eqref{eq3.10} and using \eqref{eq2.7}, \eqref{2.17}, \eqref{2.18}, \eqref{eq3.18}, we have,
\begin{align*}
\mathcal{S}\bigg{(}{_{k}^{C}\mathcal{D}}_{\alpha, \omega, 0^{+}}^{\gamma,\mu,\nu} f(t)\bigg{)}(u)&=\mathcal{S}\bigg{(}k{_{k}P}_{\alpha,(k-\mu), \omega, 0^{+}}^{-\gamma}\frac{d}{dt}f\bigg{)}(u)\\
&=k\mathcal{S}\Bigg{(}\bigg{(}{_{k}\varepsilon}_{\alpha,(k-\mu),\omega}^{-\gamma}*\frac{d}{dt}f\bigg{)}(t)\Bigg{)}(u)\\
&=k\bigg{(}\frac{u}{k}\bigg{)}^{\frac{(k-\mu)}{k}}\Bigg{(}1-\omega k\bigg{(}\frac{u}{k}\bigg{)}^{\frac{\alpha}{k}}\Bigg{)}^{\frac{\gamma}{k}}\frac{1}{u}\bigg{(}F(u)-f(0^{+})\bigg{)}
\end{align*}
Hence proved.\\
For absolutely continuous function $f \in AC^{1}[0, b]$, the relations \eqref{eq3.7} and \eqref{eq3.16} can be obtained by using Sumudu transform instead of Laplace transform technique.
\section{Applications}
In this section we find the solutions of Cauchy problems using the results obtained in above section. The problems in this section are generalizations of problems discussed in \cite{RRFZ}. For $k=1$ the following theorem \eqref{thm4.1} and theorem \eqref{thm4.2} reduces to the results in \cite{RRFZ}. 
\begin{theorem}\label{thm4.1}
The solution of Cauchy problem
\begin{align}
{_{k}\mathcal{D}}_{\alpha, \omega, 0^{+}}^{\gamma, \mu, \nu}y(x)&=\lambda _{k}P_{\alpha,\mu,\omega,0^{+}}^{\delta}y(x)+f(x), \label{eq4.1}\\
\big{[}{_{k}P}_{\alpha,(1-\nu)(k-\mu),\omega, 0^{+}}^{-\gamma(1-\nu)}y(x)\big{]}_{t=0^{+}}&=K, \qquad \qquad K\geq 0 \label{eq4.2}
\end{align}
where $x\in(0,\infty), f(x)\in L^{1}[0,\infty)$; $\mu\in(0, 1)$, $\nu\in[0, 1]$, $\omega, \lambda\in\mathbb{C}$, $\alpha > 0$, $\gamma, \delta\geq 0$ is given by,
\begin{align}
y(x)&=K\sum_{n=0}^{\infty}\lambda^{n}x^{\frac{\nu(k-\mu)+\mu(1+2n)}{k}-1}E_{k,\alpha,\nu(k-\mu)+\mu(1+2n)}^{n(\delta+\gamma)-\gamma(\nu-1)}(\omega x^{\frac{\alpha}{k}})\nonumber\\
&\quad+\sum_{n=0}^{\infty}\lambda^{n} {_{k}P}_{k,\alpha,\mu(1+2n),\omega,0^{+}}^{\gamma+n(\delta+\gamma)}f(x).\label{eq4.3}
\end{align}
\end{theorem}
\textbf{Proof:}
Let $Y(u)$ and $F(u)$ denote the Laplace transform of $ y(x) $ and $ f(x) $ respectively. Now taking Laplace transform of \eqref{eq4.1}, and using \eqref{eq2.7}, \eqref{eq2.8}, \eqref{eq2.10},\eqref{eq2.11},  \eqref{eq3.11}, \eqref{eq4.2}, we have,
\begin{align*}
&(ku)^{\frac{\mu}{k}}\bigg{(}1-\omega k(ku)^{\frac{-\alpha}{k}}\bigg{)}^{\frac{\gamma}{k}}Y(u)-k(ku)^{\frac{-\nu(k-\mu)}{k}}\bigg{(}1-\omega k(ku)^{\frac{-\alpha}{k}}\bigg{)}^{\frac{\gamma\nu}{k}}K \nonumber\\
&=\lambda (ku)^{\frac{-\mu}{k}}\bigg{(}1-\omega k(ku)^{\frac{-\alpha}{k}}\bigg{)}^{\frac{-\delta}{k}}Y(u)+F(u)
\end{align*}
On simplification,
\begin{align*}
&Y(u)=\Bigg{(}\frac{Kk(ku)^{\frac{-\nu(k-\mu)}{k}}\big{(}1-\omega k(ku)^{\frac{-\alpha}{k}}\big{)}^{\frac{\gamma\nu}{k}}+F(u)}{(ku)^{\frac{\mu}{k}}\big{(}1-\omega k(ku)^{\frac{-\alpha}{k}}\big{)}^{\frac{\gamma}{k}}-\lambda (ku)^{\frac{-\mu}{k}}\big{(}1-\omega k(ku)^{\frac{-\alpha}{k}}\big{)}^{\frac{-\delta}{k}}}\Bigg{)}\\\\
&=\Bigg{[}\frac{Kk(ku)^{\frac{-\nu(k-\mu)}{k}}\big{(}1-\omega k(ku)^{\frac{-\alpha}{k}}\big{)}^{\frac{\gamma\nu}{k}}+F(u)}{(ku)^{\frac{\mu}{k}}\big{(}1-\omega k(ku)^{\frac{-\alpha}{k}}\big{)}^{\frac{\gamma}{k}}}\Bigg{]}\frac{1}{\bigg{[}1-\frac{\lambda (ku)^{\frac{-\mu}{k}}\big{(}1-\omega k(ku)^{\frac{-\alpha}{k}}\big{)}^{\frac{-\delta}{k}}}{{(ku)^{\frac{\mu}{k}}\big{(}1-\omega k(ku)^{\frac{-\alpha}{k}}\big{)}^{\frac{\gamma}{k}}}}\bigg{]}}
\end{align*}
Therefore, for $\left|\frac{\lambda (ku)^{\frac{-\mu}{k}}\big{(}1-\omega k(ku)^{\frac{-\alpha}{k}}\big{)}^{\frac{-\delta}{k}}}{{(ku)^{\frac{\mu}{k}}\big{(}1-\omega k(ku)^{\frac{-\alpha}{k}}\big{)}^{\frac{\gamma}{k}}}}\right|<1$, we get,
\begin{align}
Y(u)&=Kk\sum_{n=0}^{\infty}\lambda^{n}(ku)^{\frac{-\nu(k-\mu)-\mu(1+2n)}{k}}\big{(}1-\omega k(ku)^{\frac{-\alpha}{k}}\big{)}^{\frac{\gamma(\nu-1)-n(\delta+\gamma)}{k}}\nonumber\\
&\quad+F(u)\sum_{n=0}^{\infty}\lambda^{n}(ku)^{\frac{-\mu(1+2n)}{k}}\big{(}1-\omega k(ku)^{\frac{-\alpha}{k}}\big{)}^{\frac{-\gamma-n(\delta+\gamma)}{k}}.\label{eq4.4}
\end{align}
Now taking inverse Laplace transform on both side of \eqref{eq4.4}, we get the required solution \eqref{eq4.3}.

\begin{theorem}\label{thm4.2}
The solution of Cauchy problem
\begin{align}
{_{k}^{C}\mathcal{D}}_{\alpha, \omega, 0^{+}}^{\gamma,\mu,\nu}u(x,t)&=K\frac{\partial^{2}}{\partial x^{2}}u(x,t),\quad\quad t>0,\quad x\in R\label{eq4.5}\\
u(x,0^{+})&=g(x),\label{eq4.6}\\
\lim_{x\to\pm\infty} u(x,t)&=0,\label{eq4.7}
\end{align}
with $\mu\in(0,1)$; $\omega\in \mathbb{R}$; $K, \alpha > 0$, $\gamma \geq 0$ is given by
\begin{align}
u(x,t)=\frac{1}{2k^{2}\pi}\int_{-\infty}^{\infty}dp\,e^{ipx}\widehat{g}(p)\sum_{n=0}^{\infty}(-K)^{n} p^{2n}t^{\frac{n\mu}{k}}E_{k,\alpha,n\mu+k}^{n\gamma}(\omega t^{\frac{\alpha}{k}}).\label{eq4.8}
\end{align}
\end{theorem}
\textbf{Proof:} Let $\overline{u}(x, q)$ and $\widehat{u}(p, t)$  denote the Laplace transform and Fourier transform of $u(x, t)$ respectively. Taking Fourier transform \eqref{eq2.2} of \eqref{eq4.5} and using \eqref{eq4.7} we get,
\begin{align}
{_{k}^{C}\mathcal{D}}_{\alpha, \omega, 0^{+}}^{\gamma,\mu,\nu}\widehat{u}(p,t)=-Kp^{2}\quad\widehat{u}(p, t)\label{eq4.9}
\end{align}
Now taking Laplace transform of \eqref{eq4.9} and using \eqref{eq3.12}, \eqref{eq4.6} we get,
\begin{align*}
(kq)^{\frac{\mu}{k}}\big{(}1-\omega k(kq)^{\frac{-\alpha}{k}}\big{)}^{\frac{\gamma}{k}}\bigg{(}\overline{\widehat{u}}(p,q)-\frac{g(x)}{q}\bigg{)}=-Kp^{2}\quad\overline{\widehat{u}}(p, q)\\
\bigg{(}(kq)^{\frac{\mu}{k}}\big{(}1-\omega k(kq)^{\frac{-\alpha}{k}}\big{)}^{\frac{\gamma}{k}}+K\,p^{2}\bigg{)}q\overline{\widehat{u}}(p, q)=(kq)^{\frac{\mu}{k}}\big{(}1-\omega k(kq)^{\frac{-\alpha}{k}}\big{)}^{\frac{\gamma}{k}}\widehat{g}(p)
\end{align*}
On simplification, we have,
\begin{align*}
\overline{\widehat{u}}(p, q)=\frac{\widehat{g}(p)}{q}\bigg{(}1+\frac{K\,p^{2}}{(kq)^{\frac{\mu}{k}}\big{(}1-\omega k(kq)^{\frac{-\alpha}{k}}\big{)}^{\frac{\gamma}{k}}}\bigg{)}^{-1}
\end{align*}
Therefore, for $\left|\frac{K\,p^{2}}{(kq)^{\frac{\mu}{k}}\big{(}1-\omega k(kq)^{\frac{-\alpha}{k}}\big{)}^{\frac{\gamma}{k}}}\right|<1$, we get,
\begin{align}
\overline{\widehat{u}}(p, q)&=\frac{\widehat{g}(p)}{q}\sum_{n=0}^{\infty}(-K)^{n} p^{2n}(kq)^{\frac{-n\mu}{k}}\big{(}1-\omega k(kq)^{\frac{-\alpha}{k}}\big{)}^{\frac{-n\gamma}{k}}\nonumber\\
&=\frac{\widehat{g}(p)}{k}\sum_{n=0}^{\infty}(-K)^{n} p^{2n}(kq)^{\frac{-n\mu-k}{k}}\big{(}1-\omega k(kq)^{\frac{-\alpha}{k}}\big{)}^{\frac{-n\gamma}{k}}.\label{eq4.10}
\end{align}
Taking inverse Laplace transform of \eqref{eq4.10}, we get,
\begin{align}
\widehat{u}(x,q)=\frac{\widehat{g}(p)}{k^{2}}\sum_{n=0}^{\infty}(-K)^{n} p^{2n}t^{\frac{n\mu}{k}}E_{k,\alpha,n\mu+k}^{n\gamma}(\omega t^{\frac{\alpha}{k}}).\label{eq4.11}
\end{align}
Taking inverse Fourier transform \eqref{eq2.3} of \eqref{eq4.11}, we get required solution \eqref{eq4.8}.\\
The theorem \eqref{thm4.1} and theorem \eqref{thm4.2} can be proved using the Sumudu transform instead of Laplace transform technique to get same solutions \eqref{eq4.3} and \eqref{eq4.8} respectively.

\end{document}